\documentclass[12pt]{amsart}
\usepackage{amssymb,bbold}
\usepackage[all]{xy}

\theoremstyle{plain}
\newtheorem{theorem}{Theorem}
\newtheorem{corollary}[theorem]{Corollary}
\newtheorem{lemma}[theorem]{Lemma}

\theoremstyle{definition}
\newtheorem{definition}[theorem]{Definition}
\newtheorem{remark}[theorem]{Remark}

\numberwithin{equation}{section} % number equation within sections

\begin{document}
\baselineskip 18pt

\title[The Grothendieck property from an ordered point of view]
      {The Grothendieck property from an ordered point of view}

%\author[R. ~Sabbagh]{Raheleh Sabbagh}
\author[O.~Zabeti]{Omid Zabeti}

%\address[V.G.~Troitsky]
%  {Department of Mathematical
 % and Statistical Sciences, University of Alberta, Edmonton,
  %AB, T6G\,2G1. Canada}
%\email{troitsky@ualberta.ca}

\address[O.~Zabeti]
  {Department of Mathematics, Faculty of Mathematics, Statistics, and Computer science,
   University of Sistan and Baluchestan, Zahedan,
   P.O. Box 98135-674. Iran}
\email{o.zabeti@gmail.com}

\keywords{The Grothendieck property, the unbounded Grothendieck property, the positive Grothendieck property, the weak Grothendieck property.}
\subjclass[2010]{Primary: 18F10. Secondary: 46B42.}

\begin{abstract}
In this note, we consider several notions related to the Grothendieck property. Among them, we introduce the notion "unbounded Grothendieck property" in a Banach lattice as an unbounded version of the known Grothedieck property in the Banach space theory. Beside other results, surprisingly, we show that spaces with the unbounded Grothendieck property are exactly the reflexive Banach lattices.

\end{abstract}

\date{\today}

\maketitle

\section{motivation and preliminaries}
Let us start with some motivation. There are several known and important concepts in the category of all Banach spaces such as the Schur property, the Banach-Saks property, the Grothendieck property and so on. When we are dealing with a Banach lattice, as a special case of Banach spaces, the order structure comes to the mind as a fruitful tool. This structure enables us to consider different types of notions related to the property using both topological and order structures. Recall that a Banach space $X$ has the Grothendieck property if for every sequence $({x_n}')\subseteq X'$, we have ${x_n}'\xrightarrow{w^{*}}0$ implies that ${x_n}'\xrightarrow{w}0$. It is known that every reflexive space has the Grothendieck property. Moreover, in the setting of the Banach lattice theory, every $\sigma$-order complete $AM$-space with unit has this property, as well ( see \cite[Theorem 4.44]{AB}). In this paper, we consider several notions related to the Grothendieck property using order and topological structures. We investigate some relations between them. In particular, we characterize reflexive Banach lattices in terms of them. Before anything, let us recall some preliminaries. Suppose $E$ is a Banach lattice. A sequence $({x_n}')\subseteq E'$ is said to be unbounded $weak^{*}$ convergent to $x'\in E'$ ( in notation, ${x_n}'\xrightarrow{uaw^{*}} x'$) provided that $|{x_n}'-x'|\wedge u'\xrightarrow{w^{*}} 0$ for each $u'\in E_{+}^{'}$. Two elements $x,y \in E$ are called disjoint if $|x|\wedge |y|=0$. For more information about unbounded convergences as well as a comprehensive context regarding Banach lattices, see \cite{AB,Z}.
% Now, we define an unbounded version of this property as follows. Suppose $E$ is a Banach lattice. $E$ is said to have the {unbounded Grothendieck property} {\bf UGP, for short} if for every bounded sequence $({x_n}')\subseteq E'$, ${x_n}'\xrightarrow{uaw^{*}}0$ implies that ${x_n}'\xrightarrow{w}0$. It is easy to see that every reflexive space possesses the unbounded Grothendieck property, as well. In spite of \cite[Theorem 4.44]{AB}, $\ell_{\infty}$ does not have the unbounded Grothendieck property; we have the following elementary lemma.
\section{main results}
First, we consider the following definition.
\begin{definition}
Suppose $E$ is a Banach lattice. $E$ is said to have
\begin{itemize}
		\item[\em (i)] {The Grothendieck  property ( {\bf GP}, for short) if for every sequence $({x_n}')\subseteq E'$, ${x_n}'\xrightarrow{w^{*}}0$ implies that ${x_n}'\xrightarrow{w}0$.}
\item[\em (ii)]
	{ The positive Grothendieck  property ( {\bf PGP}, for short) if for every sequence $({x_n}')\subseteq E_{+}^{'}$, ${x_n}'\xrightarrow{w^{*}}0$ implies that ${x_n}'\xrightarrow{w}0$.}
\item[\em (iii)]
	{The weak Grothendieck  property ( {\bf WGP}, for short) if for every disjoint sequence $({x_n}')\subseteq E'$, ${x_n}'\xrightarrow{w^{*}}0$ implies that ${x_n}'\xrightarrow{w}0$.}
		\item[\em (iv)] {The disjoint Grothendieck  property ( {\bf DGP}, for short) if for every norm bounded disjoint sequence $({x_n}')\subseteq E'$, we have  ${x_n}'\xrightarrow{w}0$.}
\item[\em (v)]{The unbounded Grothendieck  property ( {\bf UGP}, for short) if for every norm bounded sequence $({x_n}')\subseteq E'$, ${x_n}'\xrightarrow{uaw^{*}}0$ implies that ${x_n}'\xrightarrow{w}0$.}
			\end{itemize}
\end{definition}
The first part is the original definition of the Grothendieck property. The second part of the definition was initially defined in \cite{W}. The third part was proposed in \cite{MF} at first. Parts $(iv)$ and $(v)$ are new definitions using unbounded convergence and disjointness. Observe that in parts $(i),(ii), (iii)$, norm boundedness of the proposed sequence is guaranteed by the $weak^{*}$-convergence but in parts $(iv), (v)$, we need assume that the proposed sequence is norm bounded.

First, we shall consider some elementary implications between them.
\begin{lemma}\label{4}
Suppose $E$ is a Banach lattice. Then, consider the following statements.
\begin{itemize}
\item[\em (i)] {$E$ possesses the {\bf UGP}.}
\item[\em (ii)] {$E$ possesses the {\bf DGP}.}
\item[\em (iii)] {$E$ possesses the {\bf WGP}.}

\end{itemize}
Then, we have
\[(i) \Rightarrow (ii) \Rightarrow (iii)\]
\end{lemma}
\begin{proof}
$(i)\Rightarrow (ii)$. Suppose $({x_n}')\subseteq E'$ is a bounded disjoint sequence. By \cite[Lemma 2]{Z}, ${x_n}'\xrightarrow{uaw} 0$ so that ${x_n}'\xrightarrow{uaw^{*}} 0$. By the assumption, ${x_n}'\xrightarrow{w} 0$, as desired.

$(ii)\Rightarrow (iii)$. It is trivial.
\end{proof}
Observe that the opposite implications proposed in Lemma \ref{4} do not hold, in general. Consider $\ell_{\infty}$; it possesses the {\bf WGP} but it fails to have {\bf DGP}; put $u_n=(0,\ldots,0,\underbrace{1}_{2n},\underbrace{1}_{2n+1},0,\ldots)$. It is easy to see that the sequence $(u_n)\subseteq {\ell_{\infty}}'$ is disjoint. Nevertheless, it is not weakly null. Furthermore, assume that $E=\ell_1$; it possesses the {\bf DGP}; suppose $({x_n}')\subseteq \ell_{\infty}$ is a disjoint sequence so that it is weakly null by \cite[Lemma 2, Theorem 7]{Z}. But it can not have the {\bf UGP} as the $uaw^{*}$-null sequence $(v_n)$ defined via $v_n=(\underbrace{0,\ldots,0}_{n-times},1,1,\ldots)$ is not weakly null.

By \cite[Theorem 4.44]{AB}, $\ell_{\infty}$ possesses the {\bf GP}; we shall show that it fails to have the {\bf UGP}.
\begin{lemma}\label{1}
 $\ell_{\infty}$ does not have the {\bf UGP}.
\end{lemma}
\begin{proof}
It is known that $\ell_1\subseteq (\ell_{\infty})'$ so that consider the standard basis sequence $(e_n)\subseteq (\ell_{\infty})'$. It is $uaw$-null by \cite[Lemma 2]{Z} so that $ua{w^{*}}$-null. But it is not weakly null in $(\ell_{\infty})'$; in this case, it should be weakly null in $\ell_1$ which is not possible.
\end{proof}
%Moreover, recall that a Banach lattice $E$ possesses the positive Grothendieck property {\bf PGP}  provided that for every positive sequence $({x_n}')\subseteq E'$, ${x_n}'\xrightarrow{w^{*}}0$ implies that ${x_n}'\xrightarrow{w}0$. Moreover, $E$ is said to have the weak Grothendieck property ({\bf WGP})if for every disjoint $w^{*}$-null sequence $({x_n}')\subseteq E'$, we have ${x_n}'\xrightarrow{w}0$. It is easy to see that the {\bf GP} implies both the {\bf PGP} and {\bf WGP} but the converse is not true, in general; the Banach lattice $c$ consists of all convergent sequences possesses the {\bf PGP} but it fails to have the {\bf GP}. Also, $\ell_1$ possesses the {\bf WGP} but it does not have the {\bf GP}, certainly.  For more information, see \cite{MF, W}. Finally,
Note that the {\bf UGP} implies the {\bf PGP}. Now, we have the following surprising fact.
\begin{theorem}\label{3}
A Banach lattice $E$ is reflexive if and only if it possesses the {\bf UGP}.
\end{theorem}
\begin{proof}
The direct implication is trivial since $E$ is reflexive if and only if so is $E'$. For the other side, suppose $E$ has the {\bf UGP} but it fails to be reflexive. By \cite[Theorem 4.71]{AB}, $E$ has a lattice copy of either $c_0$ or $\ell_1$. Observe that neither $\ell_1$ nor $c_0$ has the {\bf UGP}; for both cases, note that ${\ell_1}'={c_0}'=\ell_{\infty}$. Consider the sequence $(u_n)\subseteq \ell_{\infty}$ defined via $u_n=(0,\ldots,0,1,\ldots)$, in which zero is appeared $n$-times. It is easy to see that by the $w^*$-topology inherited by either $\ell_1$ or $c_0$, the sequence $(u_n)$ is $ua{w^*}$-null. Nevertheless, by the Dini's theorem ( \cite[Theorem 3.52]{AB}), it is not weakly null in $\ell_{\infty}$, certainly. This shows that $\ell_1$ or $c_0$ does not have the {\bf UGP}. Now, Suppose $E$ has the {\bf UGP}. Note that by  \cite[Theorem 2.4.12]{Ni}, there is a positive projection $P$ from $E$ onto $c_0$. Moreover, by \cite[Theorem 5.3.13]{Ni}, $E$ can not have the {\bf PGP}. This also show that $E$ fails to have the {\bf UGP}, as well. Furthermore, by \cite[Proposition 2.3.11]{Ni}, there exists a positive projection $P_1$ from $E$ onto $\ell_1$. Now, consider the positive operator $T=\iota o P_1$, where $\iota$ denotes the inclusion operator from $\ell_1$ into $c_0$. Again, using \cite[Theorem 5.3.13]{Ni} convinces us that $E$ fails to have the {\bf PGP} so that it does not have the {\bf UGP}. Both statements contradict our assumption. So, $E$ is reflexive.
%\cite[Page 152. Exercise 3.42]{Fab},  so that $E$ fails to have the unbounded Grothendieck property, as well, a contradiction.
\end{proof}
Furthermore, we have the following useful observation.
%\begin{corollary}
%Suppose $E$ is a Banach lattice. Then the following implications hold.
%\[{\bf UGP} \Rightarrow {\bf GP} \Rightarrow {\bf PGP}\]
%\[{\bf UGP} \Rightarrow {\bf GP} \Rightarrow {\bf WGP}\]

%\end{corollary}
%\begin{proof}
%Suppose $({x_n}')\subseteq E'$ is a disjoint $w^{*}$-null sequence. Then,$({x_n}')$  is $uaw$-null by \cite[Lemma 2]{Z} so that $uaw^{*}$-null. By the assumption, $({x_n}')$ is weakly null as desired.
%\end{proof}
\begin{theorem}\label{2}
Suppose $E$ is a $\sigma$-order complete Banach lattice whose dual space is order continuous. Then the {\bf PGP} in $E$ implies the {\bf UGP} if and only if $E$ is order continuous.
\end{theorem}
\begin{proof}
First, assume that $E$ is order continuous and possesses the {\bf PGP}.
%Observe that by \cite[Theorem 5.3.13]{Ni}, it is enough to show that every disjoint positive $uaw^{*}$-null sequence in $E_{+}^{'}$ is also weakly null.
Suppose $({x_n}')\subseteq E'$ is a positive $uaw^{*}$-null sequence. By \cite[Theorem 4.18]{AB}, for each arbitrary $\varepsilon>0$, there exists a $u'\in E_{+}^{'}$ such that for each $x\in E_{+}$, ${x_n}'(x)-({x_n}'\wedge u')(x)<\varepsilon$. By the assumption, $({x_n}'\wedge u')(x)\rightarrow 0$ so that ${x_n}'(x)\rightarrow 0$. Therefore, ${x_n}'\xrightarrow{w^{*}}0$. By the assumption, again, ${x_n}'\xrightarrow{w}0$. For the converse, suppose not. By \cite[Theorem 4.51]{AB}, $E$ possesses a lattice copy of $\ell_{\infty}$. Note that although $\ell_{\infty}$ has the {\bf PGP} but it fails to have the {\bf UGP} by Lemma \ref{1}. This implies that $E$ does not have the {\bf UGP}, as well.
\end{proof}
\begin{corollary}
Suppose $E$ is a $\sigma$-order complete Banach lattice whose dual space is order continuous and possesses the {\bf PGP}. Then, it is reflexive if and only if it is order continuous.
\end{corollary}
Note that  the {\bf GP} always implies the {\bf PGP}. So, we have the following result.
\begin{corollary}
Suppose $E$ is a $\sigma$-order complete Banach lattice whose dual space is order continuous. Then the {\bf GP} in $E$ implies the {\bf UGP} if and only if $E$ is order continuous.
\end{corollary}
By \cite[Proposition 4.9]{MF}, when $E'$ is order continuous, the {\bf WGP} implies the {\bf PGP}. Thus, the following result comes to the mind, readily.
\begin{corollary}\label{5}
Suppose $E$ is a $\sigma$-order complete Banach lattice whose dual space is order continuous. Then the {\bf WGP} in $E$ implies the {\bf UGP} if and only if $E$ is order continuous.
\end{corollary}
\begin{remark}
Note that order continuity of $E'$ is essential in Corollary \ref{5} and can not be removed; $\ell_1$ possesses the {\bf WGP} but it fails to have {\bf UGP} although it is order continuous. Furthermore, suppose a Banach lattice $E$ possesses the {\bf UGP}. By Theorem \ref{3}, it is reflexive so that $E$ has the {\bf GP}, {\bf WGP}, and {\bf PGP}.
\end{remark}
Now, we consider the following statement; a relation between the remarkable Grothendieck property and the disjoint Grothendieck one.
\begin{theorem}
A $\sigma$-order complete Banach lattice $E$ is order continuous if and only if the {\bf GP} in $E$ implies the {\bf DGP}.
\end{theorem}
\begin{proof}
Suppose $E$ is order continuous and possesses the {\bf GP}. Assume that $({x_n}')$ is a bounded disjoint sequence in $E'$. By \cite[Lemma 2]{Z}, ${x_n}'\xrightarrow{uaw}0$ so that ${x_n}'\xrightarrow{uaw^{*}}0$. By \cite[Proposition 5]{Z}, ${x_n}'\xrightarrow{w^{*}}0$. By the {\bf GP}, ${x_n}'\xrightarrow{w}0$. For the converse, assume that $E$ is not order continuous. By \cite[Theorem 4.56]{AB}, $\ell_{\infty}$ is lattice embeddable in $E$. By \cite[Theorem 4.44]{AB}, $\ell_{\infty}$ has the {\bf GP}; nevertheless, it does not possess the {\bf DGP}; consider the standard basis sequence $(e_n)\subseteq {\ell_{\infty}}^{*}$. This implies that $E$ also can not have the {\bf DGP}, as well.
\end{proof}
%Furthermore, we have the following result regarding the relations between the positive Grothendieck property and the unbounded Grothendieck one.
%\begin{theorem}
%A Banach lattice $E$ with the {\bf UGP} possesses the {\bf PGP} if and only if $E'$ is order continuous.
%\end{theorem}
%\begin{proof}
%Suppose $E'$ is order continuous. By \cite[Theorem 5.3.13]{Ni}, it is enough to show that every disjoint positive $w^{*}$-null sequence in $E_{+}^{'}$ is also weakly null. Suppose $({x_n}')\subseteq E_{+}^{'}$ is a disjoint $w^{*}$-null sequence. By \cite[Lemma 2]{Z}, it is $uaw$-null so that $uaw^{*}$-null. By the assumption, ${x_n}'\xrightarrow{w}0$. For the other implication, suppose not; that is $E'$ is not order continuous. By \cite[Theorem 4.69]{AB}, $\ell_1$ is lattice embeddable in $E$. It is known that $\ell_1$ does not have the {\bf PGP} so that $E$ can not have the {\bf PGP}, as well.
%\end{proof}

\end{document}